# LQR Based Improved Discrete PID Controller Design via Optimum Selection of Weighting Matrices Using Fractional Order Integral Performance Index

Saptarshi Das[1,2], Indranil Pan[2], Kaushik Halder[2,3], Shantanu Das[4], and Amitava Gupta[1,2]

1) School of Nuclear Studies & Applications (SNSA), Jadavpur University, Salt Lake Campus, LB-8, Sector 3, Kolkata-700098, India. (E-mail: saptarshi@pe.jusl.ac.in, amitg@pe.jusl.ac.in).

2) Department of Power Engineering, Jadavpur University, Salt Lake Campus, LB-8, Sector 3, Kolkata-700098, India. (E-mail: indranil.jj@student.iitd.ac.in, indranil@pe.jusl.ac.in)

3) Department of Electronics and Instrumentation Engineering, National Institute of Science & Technology, Palur Hills, Berhampur-761008, Orissa, India. (E-mail: bubun85@gmail.com)

4) Reactor Control Division, Bhabha Atomic Research Centre, Mumbai-400085, India. (E-mail: shantanu@barc.gov.in)

**Abstract:**
The continuous and discrete time Linear Quadratic Regulator (LQR) theory has been used in this paper for the design of optimal analog and discrete PID controllers respectively. The PID controller gains are formulated as the optimal state-feedback gains, corresponding to the standard quadratic cost function involving the state variables and the controller effort. A real coded Genetic Algorithm (GA) has been used next to optimally find out the weighting matrices, associated with the respective optimal state-feedback regulator design while minimizing another time domain integral performance index, comprising of a weighted sum of Integral of Time multiplied Squared Error (ITSE) and the controller effort. The proposed methodology is extended for a new kind of fractional order (FO) integral performance indices. The impact of fractional order (as any arbitrary real order) cost function on the LQR tuned PID control loops is highlighted in the present work, along with the achievable cost of control. Guidelines for the choice of integral order of the performance index are given depending on the characteristics of the process, to be controlled.

## 1. Introduction

Classical optimal control theory has evolved over decades to formulate the well known Linear Quadratic Regulators which minimizes the excursion in state trajectories of a system while requiring minimum controller effort [1]. This typical behaviour of LQR has motivated control designers to use it for the tuning of PID controllers [2]-[3]. PID controllers are most common in process industries due to its simplicity, ease of implementation and robustness. Using the Lyapunov's method, the optimal quadratic regulator design problem reduces to the Algebraic Riccati Equation (ARE) which is solved to calculate the state feedback gains for a chosen set of weighting matrices. These weighting matrices regulate the penalties on the deviation in the trajectories of the state variables ($x$) and control signal ($u$). Indeed, with an arbitrary choice of weighting matrices, the classical state-feedback optimal regulators seldom show good set-point tracking performance due to the absence of integral term unlike the PID controllers. Thus, combining the tuning philosophy of PID controllers with the concept of LQR allows the designer to enjoy both optimal set-point tracking and optimal cost of control within the same design framework.

Optimal control theory has been extended to tune PID controllers in few recent literatures. In Choi and Chung [4], an inverse optimal PID controller is designed considering the error and its integro-differential as the state variables, similar to the approach, presented in this paper. In Arruda *et al.* [5], a custom cost function has been minimized with GA to design multi-loop PID controllers as the weighted sum of ITSE and variance of the manipulated variable and controlled variable. PID controller tuning with state-space approach using the error and its first and second order derivative has been investigated in [6]-[7]. The method proposed LQR-PID of He *et al.* [2]-[3] has been extended for first and second order systems with zeros in the process model in Ghartemani *et al.* [8]. Ochi and Kondo [9] have shown that the integral type optimal servo for second order system can be reduced to a LQR problem and an optimal I-PD controller can be designed with this technique. Several classical optimal and robust control approaches of PID controller can be cast into a Linear Matrix Inequality (LMI) problem as in Ge *et al.* [10] which consider the controlled variable, its rate and integral of error as the state variables.

Genetic algorithm and other stochastic global optimization techniques have also been employed for various optimal control problems. Wang *et al.* [11] used GA to optimally find out the weighting matrices of LQR i.e. $Q$ and $R$ with a specified structure. The concept of GA based optimum selection of weighting matrices has been extended for LQR as well as pole placement problems in Poodeh *et al.* [12]. GA based optimal time domain [13] and frequency domain loop-shaping [14] based PID tuning problems are also popular in the contemporary research community. The mixed $H_2/H_\infty$ optimal PID controller tuning of Chen *et al.* [14] has been improved with GA as a single objective disturbance rejection PID controller in Krohling and Rey [15] and as multi-objective loop-shaping based design in Lin *et al.* [16]. A wide class of standard optimal control problems has been solved using evolutionary and swarm intelligence based global optimization techniques in Ghosh *et al.* [17], [18].

Fractional order systems and controllers are becoming increasingly popular in the automation and process control community. A state of the art survey on the design and application of fractional order system and controllers can be found in [19]. For optimum set-point tracking control of PID/FOPID controllers, time domain performance index optimization based tuning techniques are more popular and have been applied in Cao *et al.* [20], Das *et al.* [21] and Pan *et al.* [22], [23]. The impact of choosing the weighting matrices of LQR are discussed by Saif [24] in a detailed manner. The present methodology selects the weighting matrices for the quadratic regulator design similar to that in [11], [12], using Genetic Algorithm while minimizing a suitable time domain performance index. Then a new arbitrary (fractional) order integral performance index has been used as the objective function of GA, as suggested by Romero *et al.* [25] for signal processing applications. The impact of these new FO integral indices based PID design on the closed loop control performance as well as the corresponding optimality of the quadratic regulators are also highlighted in the present work. An analog PID controller and its discretized form a digital PID both have been tuned with the proposed optimum weight selection based corresponding continuous and discrete time LQR techniques for second order systems with very low and high damping as two illustrative examples.

The rest of the paper is organized as follows. Section 2 discusses about the theoretical framework for LQR based optimal analog and digital PID controller design. Section 3 proposes the GA based optimum weight selection methodology for LQR tuning of PID controllers. Section 4 validates the proposed argument with two classes of second order systems as two illustrative examples. The paper ends with the conclusion as section 5, followed by the references.

## 2. Formulation of LQR Based Optimal PID Controller for Second Order Systems
### 2.1. Tuning of PID Controllers as Continuous Time Linear Quadratic Regulators

He *et al.* [2]-[3] has given a formulation for tuning over-damped or critically-damped second order systems having two real open loop process poles. The concept has been extended in this sub-section for lightly damped processes as well. Also, in [2], it has been suggested that one of the real poles needs to be cancelled out by placing one of the controller zeros at the same position on the negative real axis of complex s-plane. Thus the second order plant to be controlled with a PID controller can be reduced to a first order process to be controlled by a PI controller. Indeed, this approach of He *et al.* [2] does not hold for lightly damped processes having oscillatory open loop dynamics as such reduction in not possible in this case. With the approach of optimal PID tuning for second order processes in [2], also the provision of simultaneously and optimally finding the three parameters of a PID controller (i.e. $K_p, K_i, K_d$) is lost that has been addressed in this paper. The present approach assumes the error, its rate and integral as the state variables and designs the optimal state-feedback controller gains as the PID controller parameters (Fig. 1).

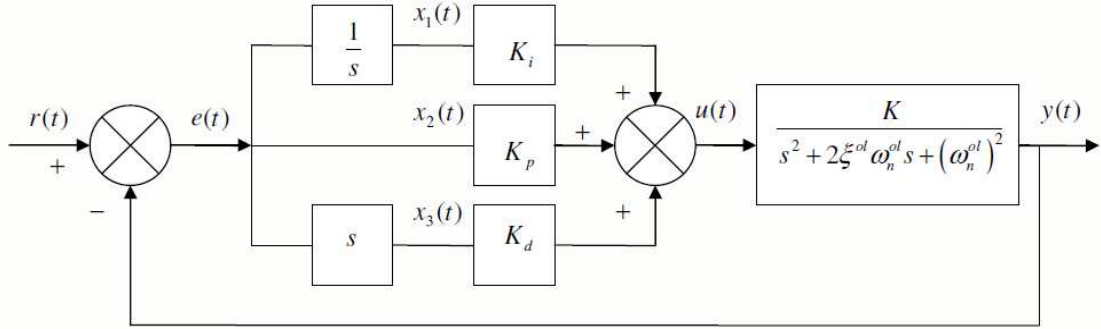

Fig.1. LQR Formulation of PID controller for second order processes.

In Fig. 1, a PID controller in parallel form (with proportional, integral and derivative gains as $K_p, K_i, K_d$) has been considered to control a second order system with known open loop damping ratio and natural frequency i.e. $\xi^{ol}, \omega_n^{ol}$ respectively. If the feedback control system is excited with an external input $r(t)$ to get a control signal $u(t)$ and output $y(t)$, then let us define the state variables as:

$$x_1 = \int e(t)dt, \quad x_2 = e(t), \quad x_3 = \frac{de(t)}{dt} \tag{1}$$

From the block diagram presented in Fig. 1, it is clear that

$$\frac{Y(s)}{U(s)} = \frac{K}{s^2 + 2\xi^{ol}\omega_n^{ol} s + \left(\omega_n^{ol}\right)^2} = \frac{-E(s)}{U(s)} \tag{2}$$

In the case of feedback design, the external set-point does not affect the controller design i.e. $r(t) = 0$. In (2), the relation $y(t) = -e(t)$ is valid for standard regulator problem as in He *et al.* [2], when the formulation is dependent on the set point. Thus, equation (2) turns out to be

$$\left[s^2 + 2\xi^{ol}\omega_n^{ol} s + \left(\omega_n^{ol}\right)^2\right] E(s) = -KU(s) \tag{3}$$

$$\Rightarrow \ddot{e} + 2\xi^{ol}\omega_n^{ol} \dot{e} + \left(\omega_n^{ol}\right)^2 e = -Ku \tag{4}$$

Using (1), equation (4) can be re-written as:

$$\dot{x}_3 + 2\xi^{ol}\omega_n^{ol} x_3 + \left(\omega_n^{ol}\right)^2 x_2 = -Ku \tag{5}$$

Using (1) and (5) the state space formulation becomes:

$$\begin{bmatrix} \dot{x}_1 \\ \dot{x}_2 \\ \dot{x}_3 \end{bmatrix} = \begin{bmatrix} 0 & 1 & 0 \\ 0 & 0 & 1 \\ 0 & -\left(\omega_n^{ol}\right)^2 & -2\xi^{ol}\omega_n^{ol} \end{bmatrix} \begin{bmatrix} x_1 \\ x_2 \\ x_3 \end{bmatrix} + \begin{bmatrix} 0 \\ 0 \\ -K \end{bmatrix} u \tag{6}$$

Comparing (6) with the standard state-space representation i.e.
$$\dot{x}(t) = Ax(t) + Bu(t) \tag{7}$$
we get the system matrices as:

$$A = \begin{bmatrix} 0 & 1 & 0 \\ 0 & 0 & 1 \\ 0 & -\left(\omega_n^{ol}\right)^2 & -2\xi^{ol}\omega_n^{ol} \end{bmatrix}, \quad B = \begin{bmatrix} 0 \\ 0 \\ -K \end{bmatrix} \tag{8}$$

In order to have a LQR formulation with the system (7), the following quadratic cost function ($J$) is minimized

$$J = \int_0^\infty \left[ x^T(t)Qx(t) + u^T(t)Ru(t) \right] dt \tag{9}$$

It has been shown in [26] that minimization of cost function (9) gives the state feedback control law as:
$$u(t) = -R^{-1}B^T Px(t) = -Fx(t) \tag{10}$$
where, $P$ is the symmetric positive definite solution of the Continuous Algebraic Riccati Equation (CARE) given by (11)
$$A^T P + PA - PBR^{-1}B^T P + Q = 0 \tag{11}$$

Here, the weighting matrix $Q$ is symmetric positive semi-definite and the weighting factor $R$ is a positive number. It is a common practice in optimal control to design regulators with varying $Q$, while keeping $R$ fixed [24] and generally they are designed with user specified closed loop performance specifications [2]. Here, Solution of the CARE (11) the $P$ matrix is generally symmetric and weighting matrix $Q$ is chosen to be a diagonal one whose elements needs to be set by the designer:

$$P = \begin{bmatrix} P_{11} & P_{12} & P_{13} \\ P_{12} & P_{22} & P_{23} \\ P_{13} & P_{23} & P_{33} \end{bmatrix}, \quad Q = \begin{bmatrix} Q_1 & 0 & 0 \\ 0 & Q_2 & 0 \\ 0 & 0 & Q_3 \end{bmatrix} \quad (12)$$

To impose high penalty on a specific state variables the elements of the matrix $Q$ can be chosen intuitively. Also, weighting factor $R$ regulates the penalty on the control signal to prevent actuator saturation. In most cases, these weights $Q, R$ are chosen with designer's expertise from the understanding of the process states. The formulation can be easily made as optimal with the use of some global optimization algorithm which will search for the weighting matrices. These optimized weights will produce the optimal regulator which will also produce an optimal time domain performance which is the main focus of this paper.

If it is now considered that the unique solution of the CARE (11) be $P$, the state feedback gain matrix becomes (13), corresponding to the optimal control signal involving the states as the loop error and its integro-differentials (1).

$$F = R^{-1}B^T P = R^{-1}\begin{bmatrix} 0 & 0 & -K \end{bmatrix} \begin{bmatrix} P_{11} & P_{12} & P_{13} \\ P_{12} & P_{22} & P_{23} \\ P_{13} & P_{23} & P_{33} \end{bmatrix}$$

$$= -R^{-1}K \begin{bmatrix} P_{13} & P_{23} & P_{33} \end{bmatrix} \quad (13)$$

$$= -\begin{bmatrix} K_i & K_p & K_d \end{bmatrix}$$

Using (10), the corresponding expression for the state feedback control signal can be derived as the output of PID controller:

$$u(t) = -Fx(t)$$

$$= -\begin{bmatrix} -K_i & -K_p & -K_d \end{bmatrix} \begin{bmatrix} x_1(t) \\ x_2(t) \\ x_3(t) \end{bmatrix} \quad (14)$$

$$= K_i \int e(t)dt + K_p e(t) + K_d \frac{de(t)}{dt}$$

The above formulation clearly shows that with judicious choice of weighting matrices $\{Q, R\}$ a PID controller can easily be tuned which preserves the achievable performance of an LQR i.e. minimum deviation in the state trajectories with minimum controller effort. The GA based choice of $Q, R$ and its advantages are discussed in the next section.

*2.2. Discrete Time Quadratic Regulator Theory Applied to Optimal Digital PID Controller Design*

It is well known that discrete time realization of PID controllers are now more preferred than their continuous time counterpart [13], [27] since the gains of a digital PID controller can be changed, switched or scheduled online so as to control complicated time varying processes over the fixed gain, lossy analog realization. In discrete time, the PID controller structure takes the following form and its performance is heavily dependent on the sampling time ($T_s$).

$$C(z) = K_p + \frac{K_i}{(1-z^{-1})} + K_d(1-z^{-1}) \qquad (15)$$

Here, the discrete time poles and zeros of the controller (15) are related with the continuous time poles and zeros by the following relation commonly known as bilinear transformation:

$$z = e^{sT_s} \Rightarrow s = \frac{2}{T_s}\left(\frac{1-z^{-1}}{1+z^{-1}}\right) \qquad (16)$$

Now, in order to do design a discrete PID controller controlling the same second order plant, the augmented system matrices (8) needs to be discretized with the specified sampling-time ($T_s$). Also, the digital PID controller should be designed with the discrete version of the optimal regulator theory to achieve optimal performance with respect to a quadratic cost function. In fact, the continuous time LQR based analog PID may not remain optimal upon discretization with arbitrary sampling time. For this reason it is a necessity to design optimal discrete PID controller with discrete version of the LQR technique.

The basics of discrete time optimal quadratic regulator is introduced here [27]. For the continuous time augmented system (8), the task is to design an optimal discrete time state feedback controller that minimizes the infinite horizon quadratic optimal cost

$$\hat{J} = \sum_{k=0}^{\infty}\left[x^T(k)Qx(k) + u^T(k)Ru(k)\right] \qquad (17)$$

Minimization of the quadratic cost given in (17) leads to the solution of the Discrete Algebraic Riccati Equation (DARE) given by (18)

$$P = Q + G^T PG - G^T PH\left(R + H^T PH\right)^{-1} H^T PG \qquad (18)$$

In (18), $Q, R$ are the positive semi-definite weighting matrices and $P$ is the positive definite solution of the discrete time Algebraic Riccati equation. The discretized system matrices can be obtained from (8) using the specified sampling time $T_s$ as:

$$\left.\begin{array}{l} G = e^{AT_s} \\ H = \left(\int_0^{T_s} e^{A\lambda}d\lambda\right)B = \left(e^{AT_s} - I\right)A^{-1}B \end{array}\right\} \qquad (19)$$

Matrix $P$ in (18) produces the optimal discrete time state-feedback gain matrix $F$ which minimizes the discrete time quadratic cost function (17) using the following relation, similar to the continuous time treatments (13):

$$F = \left(R + H^T PH\right)^{-1} H^T PG \qquad (20)$$

Thus the optimal discrete time control law is given by:

$$\begin{aligned} u(k) &= -Fx(k) \\ &= -\left(R + H^T PH\right)^{-1} H^T PGx(k) \end{aligned} \qquad (21)$$

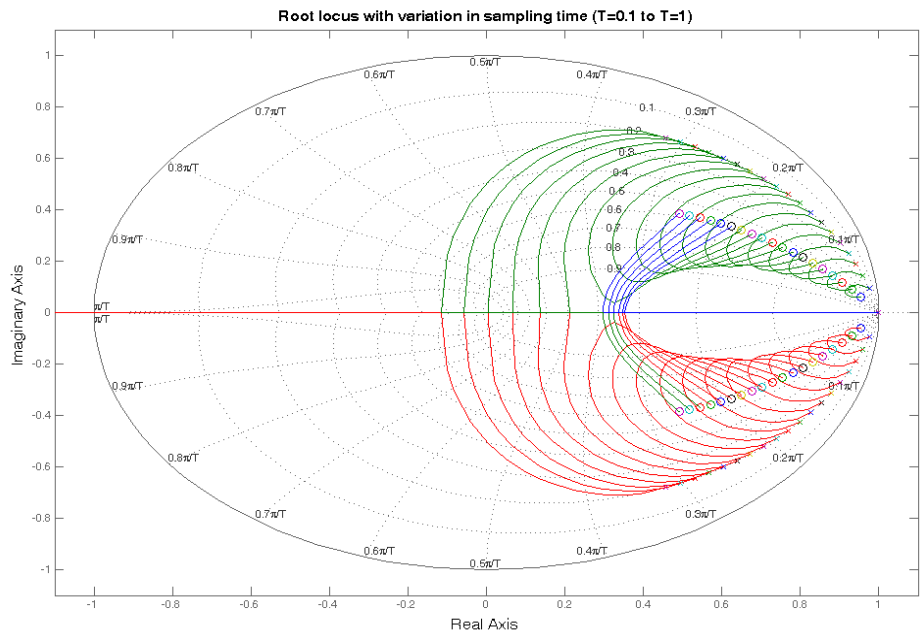

Fig. 2. Root locus of the discretized open loop system with increasing sampling time.

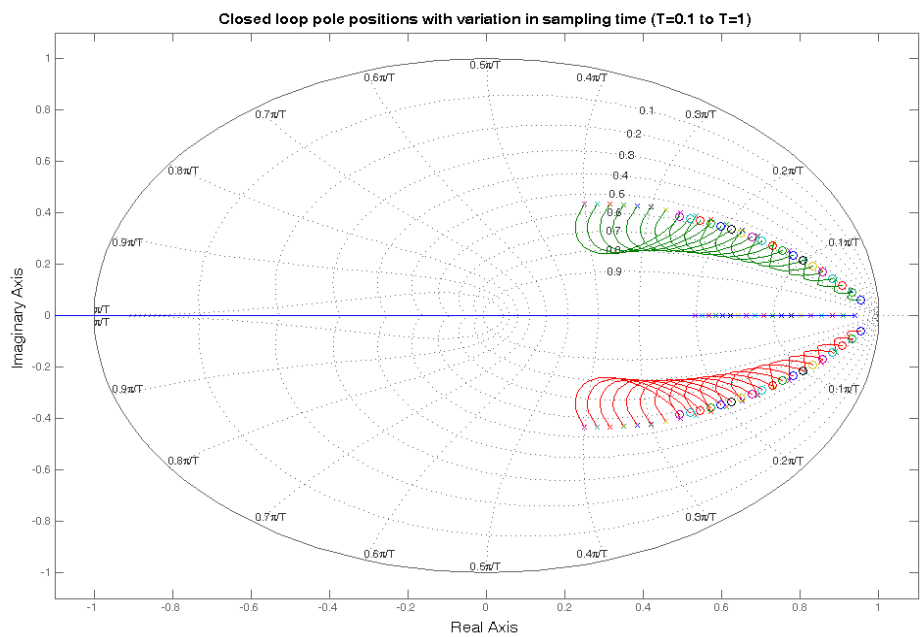

Fig. 3. Shifting of closed loop pole locations with increasing sampling time.

To demonstrate the necessity of discrete LQR based design of optimal discrete PID controller, a continuous time second order system of the structure (2) has been considered with specified gain, damping and frequency as $K=1, \xi^{ol}=0.2, \omega_n^{ol}=1$. For the continuous time regulator design with CARE (11) the weighting matrices have been considered as $Q = \begin{bmatrix} 1 & 0 & 0 \\ 0 & 1 & 0 \\ 0 & 0 & 1 \end{bmatrix}, R = 1$. Next, the control system with the resulting continuous time PID controller (13) is discretized with sampling time $T_s \in [0.1, 1]$. The open loop root locus and closed loop pole locations are shown in Fig. 2 and Fig. 3 respectively. It is clear that the dominant complex poles and PID controller zeros shift towards high frequencies thus loosing its dominant dynamic behavior. Thus, for a specific sampling time $T_s$, the optimal controller needs to be derived using the discrete version of the LQR formulation i.e. DARE given by (18).

## 3. LQR Based PID Controller Design with Optimum Selection of Weighting Matrices
### *3.1. Effect of Weighting Matrices on the Control System Performance*

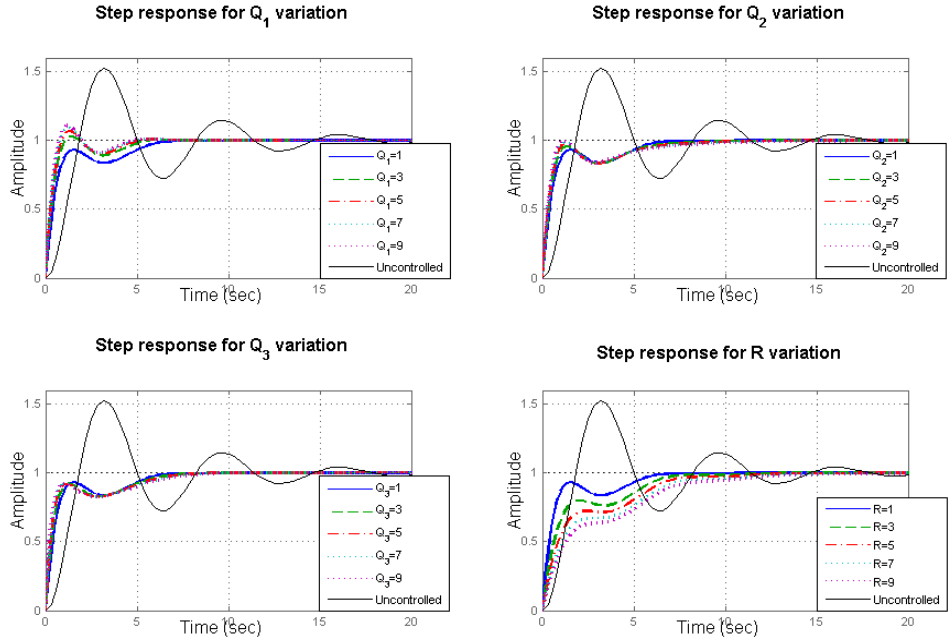

Fig. 4. Effect of weighting matrices Q and R on the time response.

Fig. 4 shows the variation in time domain performances for the example, discussed in the previous section with variation in the weighting matrices $Q, R$. With the variation in elements of matrix $Q$, the overshoot slightly increases with gradual fall in rise time. For high value of $R$, the time response becomes sluggish. Similar observations could have been found from control signal point of view. Therefore, a Genetic Algorithm based approach has been adopted in this paper for optimum choice of weighting matrices $Q, R$.

## 3.2. Optimum Selection of Weighting Matrices Using Genetic Algorithm Minimizing Time Domain Integral Performance Index

It is well known [26], the above LQR based PID controller is the most optimal for a specific choice of the weighting matrices $Q, R$. Indeed, the time domain performance is heavily affected for any arbitrary choice of the weighting matrices (Fig. 4) although the optimality in terms of the trade-off between excursion of the states trajectories and controller effort is preserved by the LQR itself. Therefore, it is logical to choose the weighting matrices in an optimum fashion with respect to another time domain performance index as these weighting matrices determine the state feedback gains (PID controller gains in this case) while indirectly monitoring the closed loop performance. Thus, a GA based stochastic optimization is formed by minimizing the cost function $\tilde{J}$ (22) as a weighted sum of ITSE and Integral Squared Controller Output (ISCO) as in [22]-[23]. This tunes the elements of the weighting matrices i.e. $\{Q_1, Q_2, Q_3, R\}$ of LQR producing time domain optimal PID.

$$\tilde{J} = \int_0^\infty \left[ w_1 \cdot t \cdot e^2(t) + w_2 \cdot u^2(t) \right] dt \tag{22}$$

Here, $w_1, w_2$ are the corresponding weights of ITSE and ISCO and are considered to be same, so as to put equal penalties on the loop error index and control signal. The rationale for using both these parameters in the objective function is to get a good time domain response and at the same time to limit the controller output to avoid actuator saturation and integral wind-up. At a first glance this might seem as a redundant repetition since the LQR methodology already gives optimal values of the controller gains with the lowest cost. However, this is actually obtained for a specified value of the weighting matrices. When $Q, R$ are varied, for each choice of weighting matrices the LQR would give an optimal gain with the lowest possible cost, but that does not necessarily imply a good time domain performance. Also, for an optimal choice of weighting matrices, the PID tuning problem becomes optimal due to the introduction of time domain performance index (22) as well as the continuous/discrete time optimal regulator (LQR) based approach (9) and (17) respectively involving the state variables.

## 3.3. Details of the Genetic Algorithm for Optimal Controller Design

Genetic Algorithm is a stochastic optimization algorithm and has been widely employed in the tuning of PID controllers, subjected to the minimization of a certain cost function like in [13]-[16], [20]-[23]. Genetic algorithm has certain advantages over the classical gradient based optimization algorithms since they are stochastic in nature and are less susceptible to get trapped in the local minima within the search space. Initially a random population of genes (which is essentially a vector comprising of the decision variables) is chosen from the search space. They undergo reproduction, crossover and mutation to yield individuals with better fitness (lower $\tilde{J}$ value in this case). A scaling function is converts the raw fitness scores in a form that is suitable for the selection function. Rank fitness scaling is used which scales the raw scores on the basis of its position in the sorted score list. This removes the effect of the spread of the raw scores. The individuals with higher fitness values have more probability of creating their copies in the next generation. This is termed as reproduction. Two parent individuals can do information interchange in a probabilistic fashion to create a child in the next generation. This process is known as crossover. In mutation a small part of the parent gene is randomly changed to yield a child. Another factor called the elite count is used which dictates the number of fittest individuals that would definitely go to the next generation. This is generally kept small with respect to the overall population so that the dominance of fitter individuals at the beginning of the simulation is reduced and premature convergence is avoided. In this case, the population size is considered to be 20 and elite count as 2.

The crossover and mutation fractions of the whole population determine the number of individuals, other than the elite, who evolve through crossover and the number which evolve through mutation respectively. This can be pre-specified by the user. The choice of these values depends on the type of optimization. In the present simulation a crossover fraction of 0.8 and a mutation fraction of 0.2 are chosen which works well for a wide variety of problems [28]. A scattered crossover function is used which creates a random binary vector and selects the genes where the vector has a value of 1 from the first parent, and the genes where the vector has a value of 0 from the second parent, and combines the genes to form the child. For mutation the Gaussian function is used which adds a random number to each vector entry of an individual. This random number is taken from a Gaussian distribution centered around zero. The other parameters of GA like population size, scaling function, selection function, elite count, mutation function, crossover function, which are used in the simulation studies, are also chosen in the lines of the previous argument. Also, a high penalty is imposed when the choices of controller gains give an unstable response. The algorithm is terminated if the maximum number of iterations is reached or the change in the objective function is lower than a specified tolerance level.

The variables that constitute the search space for the PID controller are $\{Q_1, Q_2, Q_3, R\}$. The intervals of the search space for these variables are $\{Q_1, Q_2, Q_3, R\} \in [0, 100]$. The variables are encoded as real values in the algorithm. The algorithm has also been run several times to ensure that the true global minima is found in the search space and the best results having the lowest cost function (along with the corresponding decision variables) have been reported here.

### *3.4. Fractional Order Integral Performance Indices and Their Impact on the LQR Based PID Design*

Fractional calculus is a 300 year's old subject and has found wide application in many branches of engineering and science [29]-[31]. The fractional order integral of any arbitrary function $f(t)$ can be represented by the left sided Riemann-Liouville definition as:

$$_0I_t^\alpha f(t) = \frac{1}{\Gamma(\alpha)} \int_0^t f(\tau)(t-\tau)^{\alpha-1} d\tau, \quad t \geq 0 \qquad (23)$$

The time variable $\tau$ in equation (23) can be replaced by a scale transformation $g_t(\tau)$, i.e. $\tau \to g_t(\tau)$. Taking

$$g_t(\tau) = \frac{1}{\Gamma(\alpha+1)}\{t^\alpha - (t-\tau)^\alpha\} \tag{24}$$

we have, $dg_t(\tau) = \frac{(t-\tau)^{\alpha-1}}{\Gamma(\alpha)}$ (25)

Hence, (23) can be written in the form

$$_0I_t^\alpha f(t) = \int_0^t f(\tau)dg_t(\tau) \tag{26}$$

Keeping, $t$ fixed, let us consider a 3D curve in the space $(\tau, g, f)$ given by the following expression $C_t : (\tau, g_t(\tau), f(\tau))$, $0 \leq \tau \leq t$. Along the curve $C_t$ if a "fence" is built perpendicular to the plane $(\tau, g)$ of varying height $f(\tau)$ as in Fig. 3, then the shadows cast on the walls by the fence may be interpreted as follows.

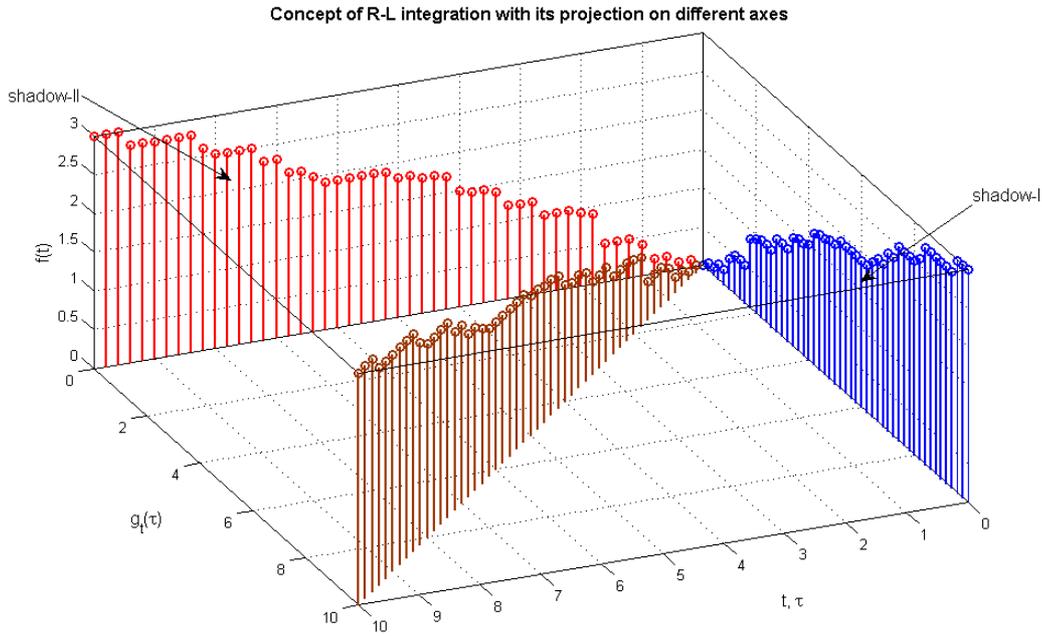

Fig. 5. Geometric interpretation of fractional order integration.

a) The area of the projection of the fence on the plane $(\tau, f)$ is given by

$$I_0^1 f(t) = \int_0^t f(\tau)d\tau \tag{27}$$

b) The area of the projection of the fence on the plane $(g, f)$ is given by

$$I_0^\alpha f(t) = \int_0^t f(\tau)dg(\tau) \tag{28}$$

For $\alpha = 1$, equation (24) reduces to $g_t(\tau) = \tau$ and hence both the shadows are equal. Thus the integer order definite integration is a special case of the left sided R-L fractional integration even from a geometric perspective. When $t$ is changing, the fence changes in length and shape. The corresponding changes in the shadow on the walls $(g, f)$ (as shown in Fig. 6) due to the change of the fence with time give a dynamical geometric interpretation of the fractional integral given by equation (23) as a function of the variable $t$.

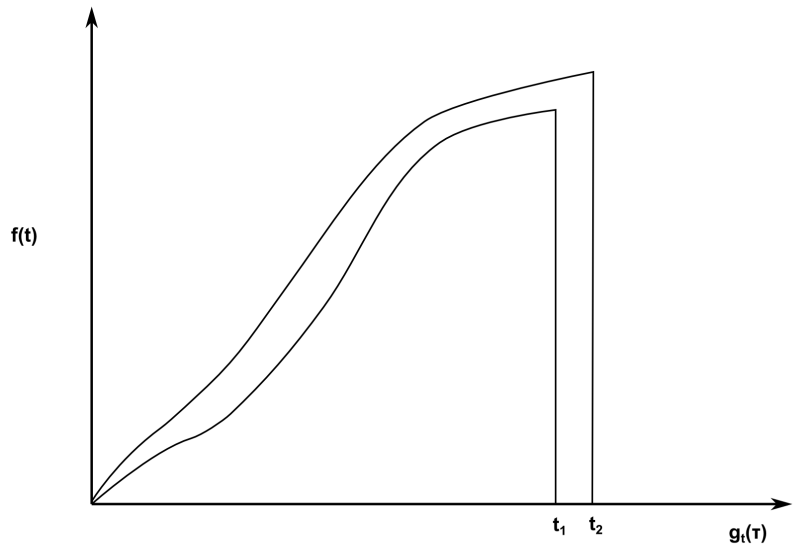

Fig. 6. Changing shadow on the wall as t changes.

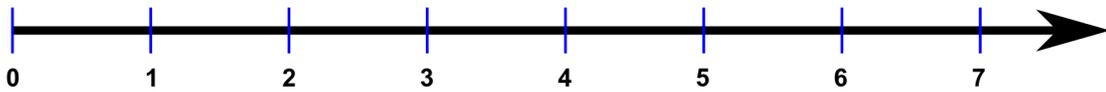

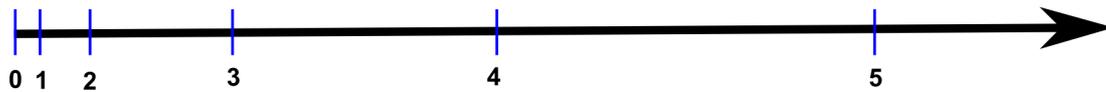

Fig. 7. The concept of homogeneous and heterogeneous time.

The physical interpretation of $I_0^\alpha f$ can be obtained by introducing the notion of a transformed time scale where the time does not flow homogeneously. Considering $\tau$ as the time, the third dimension $g_t(\tau)$ added to the pair $(\tau, f(\tau))$ can be considered as some kind of a transformed time scale. Thus the concept of two kinds of time arises as represented in Fig 7.
  a) The mathematical time $\tau$ which is assumed to be homogeneous and equably flowing.
  b) The transformed time $g(\tau)$ whose notion can be understood from the following. Assuming a clock displays the time $\tau$ incorrectly and the relationship between the measured time $\tau$ and the real time $T$ (i.e. the correct or transformed time) is given by $T = g(\tau)$.

Hence, when a real time interval of $dT = dg(\tau)$ elapses, the time interval measured using the notion of mathematical time is $d\tau$. Thus if $v(\tau)$ is the measured velocity of a body, then the wrong value of distance covered is given by the integral

$$I_0^1 v(t) = \int_0^t v(\tau) d\tau \tag{29}$$

whereas the real or actual distance passed is given by

$$I_0^\alpha v(t) = \int_0^t v(\tau) dg(\tau) \tag{30}$$

Gutierrez *et al.* [30] and Podlubny [32] have given the geometric illustration of fractional order differentiation and integration in a lucid manner. Now, a new cost function is proposed here with the generalization of the order of integration to be any arbitrary number ($\Lambda$) [25].

$$\tilde{J}^* = \frac{d^{-\Lambda}}{dt^{-\Lambda}} \left( w_1 \cdot t \cdot e^2(t) + w_2 \cdot u^2(t) \right) \tag{31}$$

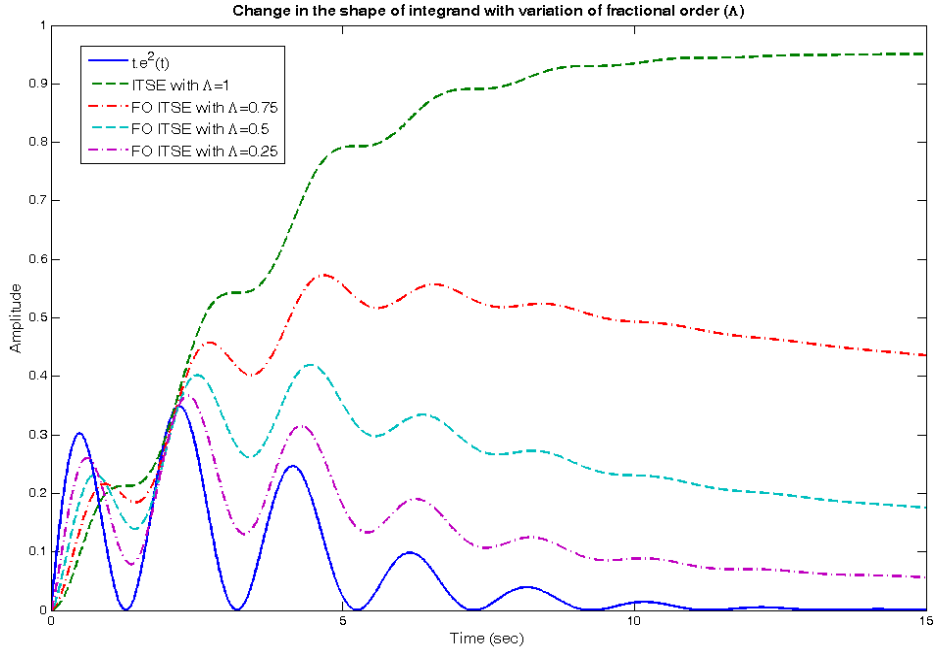

Fig. 8. Changing shape of integrand (error index) with variation in FO of integration.

This new concept of monitoring the order of cost function as an extra design tool highlights the intricacies in the functional form of the cost function perhaps in a new resolution. It is therefore logical to get superior and inferior control performance depending on the order of cost function ($\Lambda$) for higher and lower than unity. In Chapter 5 of [29], it has been illustrated in a detailed manner that for the FO integrals like that represented by (24), the integrand changes its shape with time unlike evaluation of the area, under a constant curve for integer order integrals. The FO integral represents an area under a transformed function which is dependent on the time step and fractional order with changing limits [29].

Now, to show the flexibility of fractional order integral in the error index in controller design an oscillatory system is considered with $K=1, \xi^{ol}=0.2, \omega_n^{ol}=1$ with a badly tuned PID controller setting $K_p=2, K_i=2, K_d=1$ (as a guess solution in the real coded GA based optimization process and the time evolution of ITSE has been shown in Fig. 8. Fig. 8 also shows that the ITSE which is based on a first order integration approaches towards a steady value monotonically as time increases whereas with a low order of fractional integration ($\Lambda$) the integration no longer remains monotonic function which justifies the fact that the integrand in (31) is changing its shape over time [29], [32]. A few interesting observations can be inferred from Fig. 8. The $t \cdot e^2(t)$ curve is always positive but has an oscillatory behavior and it tends towards zero with increase in time. This implies that the process settles down to a steady state value after a transitory period. Now from our concept of integer order integration, which interprets it simply as an area under a curve, we know that it is monotonically increasing if the area is positive as also shown in Fig. 8 for $\Lambda=1$. But due to the complicated characteristics of the fractional order integration, it can be seen that the integrand curves are not monotonically increasing and shows some kind of oscillatory behavior which is not possible with integer order integrals. Also the final value of the fractional order integrands vary with time and none of them reach a steady state value like the integer order integral. Thus taking the final value of the integral in a finite time horizon for the design of the controller or comparing two designs based on the precept of final values may not be appropriate.

For numerical simulation, Riemann-Liouville definition cannot be applied directly. Hence, band limited realization of fractional order differ-integrators needs to be adopted. In the present simulation study, each fractional order element has been rationalized with Oustaloup's recursive filter [33], given by the following equation (32)-(22). If it be assumed that the expected fitting range or frequency range of controller operation is $(\omega_b, \omega_h)$, then the higher order filter which approximates the FO element $s^\gamma$ can be written as:

$$G_f(s) = s^\gamma = K \prod_{k=-N}^{N} \frac{s+\omega_k'}{s+\omega_k} \qquad (32)$$

where the poles, zeros, and gain of the filter can be evaluated as:

$$\omega_k = \omega_b \left(\frac{\omega_h}{\omega_b}\right)^{\frac{k+N+\frac{1}{2}(1+\gamma)}{2N+1}}, \omega_k' = \omega_b \left(\frac{\omega_h}{\omega_b}\right)^{\frac{k+N+\frac{1}{2}(1-\gamma)}{2N+1}}, K = \omega_h^\gamma \qquad (33)$$

In (32) and (34), $\gamma$ is the minimal part of the FO order of the differ-integration and $(2N+1)$ is the order of the filter. Present study considers a $5^{th}$ order Oustaloup's rational approximation for the FO elements within the frequency range $\omega \in \{10^{-2}, 10^2\}$ rad/sec.

## 4. Simulations and Results

To show the effectiveness of the proposed methodology a heavily oscillatory and a sluggish system has been considered of the structure (2) with parameters $K=1, \xi^{ol}=\{0.2, 5\}, \omega_n^{ol}=1$ respectively [34], excluding the time delay. GA based selection of weighting matrices yields the PID controller gains as the optimal state-feedback gains for the continuous and discrete version of the LQR formulation and has been reported in Table I and II respectively along with the minima of the FO cost function ($J_{min}$). It is observed from the Tables I-II that the obtained minima of the cost function or $J_{min}$ in each cases of integer or fractional order integral is less for the CARE based design than the DARE based one since the tracking becomes slightly worse for consideration of the sampled data system which is much practical from implementation point of view. Here, all simulations have been run for a finite time horizon of 100 seconds.

Table I: Optimum PID Controller Parameters with CARE

| Process | Fractional Order of Performance Index (Λ) | $J_{min}$ | $K_p$ | $K_i$ | $K_d$ |
|---|---|---|---|---|---|
| Oscillatory | 0.5 | 14.163 | 1.366073 | 0.260574 | 1.776595 |
| | 1.0 | 96.855 | 2.291051 | 0.234846 | 3.793601 |
| | 1.5 | 842.254 | 1.802627 | 0.246183 | 2.408207 |
| Sluggish | 0.5 | 16.366 | 1.641867 | 0.172061 | 0.168591 |
| | 1.0 | 118.156 | 1.732432 | 0.16095 | 0.173052 |
| | 1.5 | 999.961 | 1.892469 | 0.252889 | 0.198146 |

Table II: Optimum PID Controller Parameters with DARE

| Process | Fractional Order of Performance Index (Λ) | $J_{min}$ | $K_p$ | $K_i$ | $K_d$ |
|---|---|---|---|---|---|
| Oscillatory | 0.5 | 14.85861 | 0.098233 | 0.171755 | 0.198904 |
| | 1.0 | 103.9139 | 0.097949 | 0.170926 | 0.198634 |
| | 1.5 | 940.4062 | 0.101853 | 0.175119 | 0.204367 |
| Sluggish | 0.5 | 16.65402 | 1.351018 | 0.16446 | 0.134639 |
| | 1.0 | 120.8922 | 1.527714 | 0.14705 | 0.151819 |
| | 1.5 | 1173.427 | 1.389176 | 0.150055 | 0.138047 |

Since the performance index itself is different for the designs (due to change of the fractional order Λ), hence numerically comparing the values of $J_{min}$ in Tables 1 and 2 is not justified. Thus some other criteria must be employed to determine which of these performance indices actually gives a better design. This has been done in the later part of this section by comparing the cost of control of the designed LQR controllers, from the obtained Riccati solutions or $P$ matrices. Another observation from the Tables 1 and 2 is that, as the fractional order of the integrand increases, the numerical value of $J_{min}$ increases. Though the original function within the integral ($t \cdot e^2(t)$) might be same, integration of the function with $\Lambda < 1$, would give a lower final value of $J_{min}$ than that given by $\Lambda = 1$ as shown in Fig. 8. For $\Lambda > 1$, the values would be higher than that given by $\Lambda = 1$ due to the inherent nature of the memory effect present in the fractional integral [29]-[32].

The time responses with the optimal PID controllers to control the oscillatory process have been compared in Fig. 9 and the associated control signals in Fig. 10. It is seen from Fig. 9 that the load disturbance rejection performance of the continuous time LQR based PID controllers are better than their discrete versions since the sensitivity reduction issues have not been considered in the optimization process. However, the set-point tracking performances are satisfactory for such a lightly damped open loop process with both types of controllers. Also, the discrete LQR based PID controllers give low and smoother control action (Fig. 10) compared to the continuous time optimal PID controllers, This especially important to make the actuator jerk free.

Similar observations can be made from Fig. 11-12, showing the time response and control signals corresponding to the sluggish plant. It is evident from Fig. 11 that the discrete LQR based PID controllers give dead-beat set-point tracking while the load disturbance performances are also comparable with the continuous time counterparts unlike that with the oscillatory system. Also, incorporating digital PID controller designed with discrete LQR one can achieve small and smooth control signals (Fig. 12) which helps to reduce size of the actuator and the associated costs.

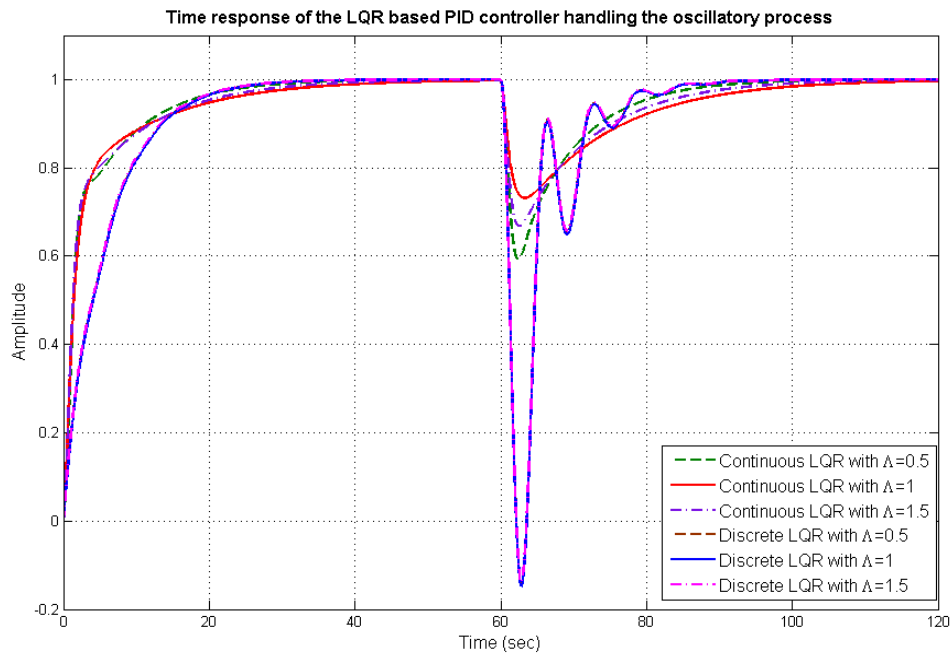

Fig. 9. Unit Set-point and load disturbance response for the oscillatory plant.

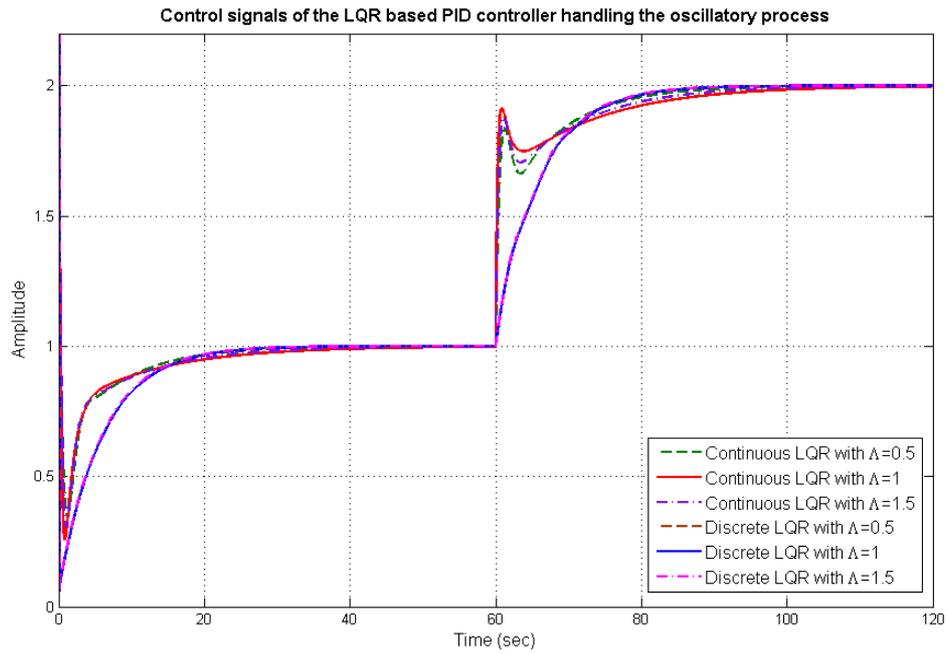

Fig. 10. Control signals of the optimal PID controllers for the oscillatory plant.

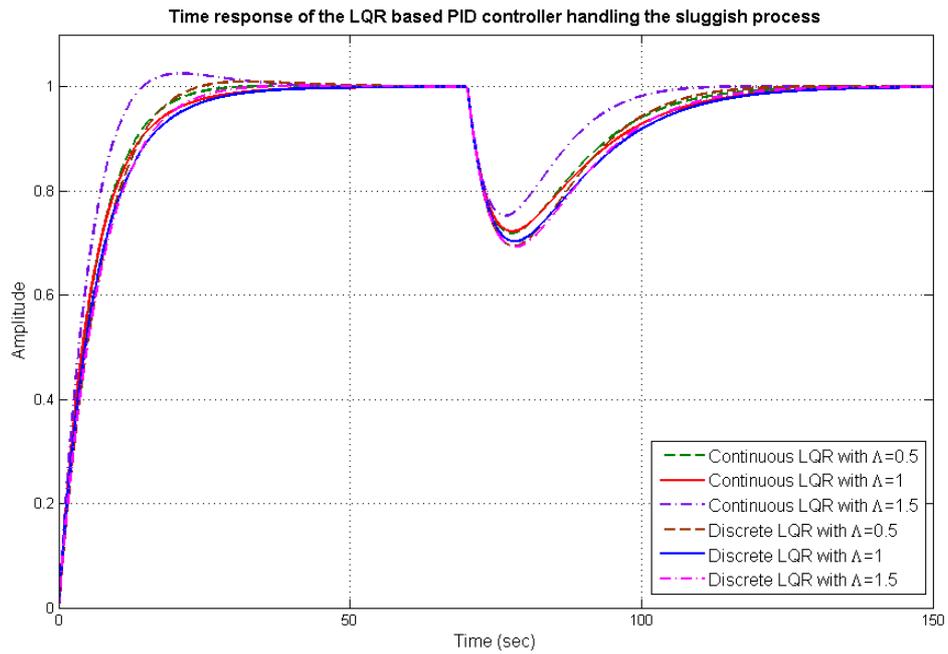

Fig. 11. Unit Set-point and load disturbance response for the sluggish plant.

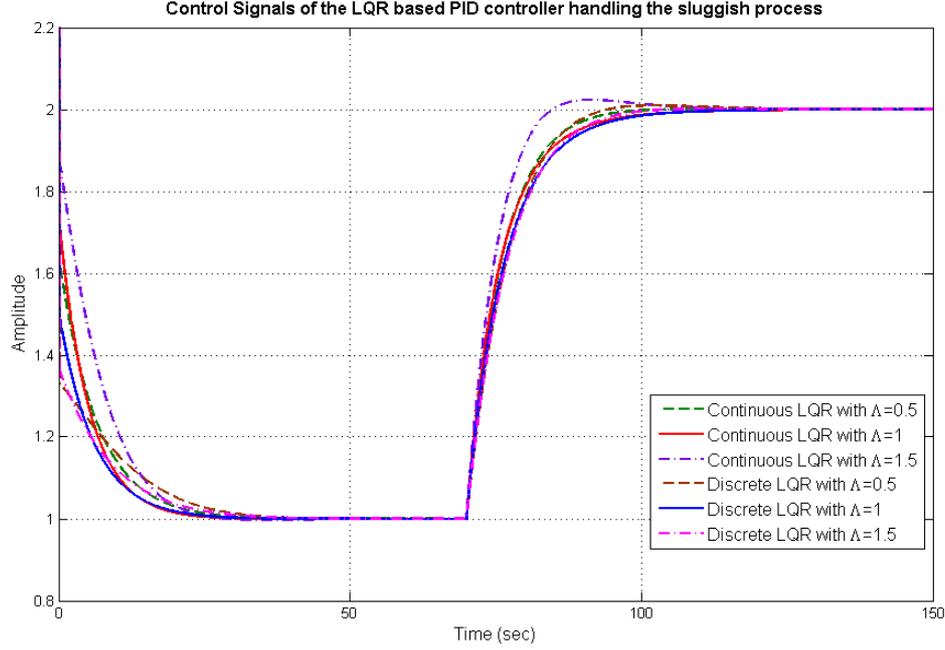

Fig. 12. Control signals of the optimal PID controllers for the sluggish plant.

From the above figures it has been seen the time responses and control signals are significantly different for the continuous and discrete time LQR based controllers but are almost same for variation with the order of fractional order integral. The rationale of incorporating the FO of the time domain integral performance index as an extra design freedom has been illustrated here. Since, the solution methodology is to search for a set of weighting matrices of LQR which produces good time response in terms of low error index and controller effort, each controller state feedback gain (discrete PID parameters) in Table II is associated with a unique set of weighting matrices satisfying the discrete algebraic Riccati equation (18). These GA based selection of weighting matrices would produce a unique Riccati solution given be the matrix $P$ which is also a measure of the associated cost of control for the LQR. The DARE based Riccati solutions for controlling the lightly damped process has been reported in (34)-(36) for varying level of FO of cost function ($\Lambda$), representing the same cost function in different resolutions.

$$P_{\Lambda=0.5} = \begin{bmatrix} 38.3375 & 20.8896 & 34.8736 \\ 20.8896 & 17.4583 & 19.9732 \\ 34.8736 & 19.9732 & 40.3691 \end{bmatrix} \quad (34)$$

$$P_{\Lambda=1.0} = \begin{bmatrix} 56.4828 & 30.7708 & 51.3929 \\ 30.7708 & 25.9626 & 29.4921 \\ 51.3929 & 29.4921 & 59.6989 \end{bmatrix} \quad (35)$$

$$P_{\Lambda=1.5} = \begin{bmatrix} 61.3195 & 33.6054 & 55.5945 \\ 33.6054 & 28.8308 & 32.3811 \\ 55.5945 & 32.3811 & 64.8493 \end{bmatrix} \quad (36)$$

Also, Riccati solutions associated with the DARE based discrete PID controller gains have been reported in (37)-(39) for varying level of $\Lambda$.

$$P_{\Lambda=0.5} = \begin{bmatrix} 247.25 & 1065.1 & 105.0964 \\ 1065.1 & 8733.6 & 863.2543 \\ 105.0964 & 863.2543 & 86.0539 \end{bmatrix} \quad (37)$$

$$P_{\Lambda=1.0} = \begin{bmatrix} 42.8740 & 172.0613 & 16.9487 \\ 172.0613 & 1789.5 & 176.0818 \\ 16.9487 & 176.0818 & 17.4946 \end{bmatrix} \quad (38)$$

$$P_{\Lambda=1.5} = \begin{bmatrix} 40.2110 & 170.5117 & 16.8189 \\ 170.5117 & 1578.077 & 155.6976 \\ 16.8189 & 155.6976 & 15.4685 \end{bmatrix} \quad (39)$$

It is shown in [26] that the infinite time performance index (17) can be calculated from the Riccati solution ($P$ matrix) using the initial values of the state variables i.e.

$$\hat{J} = x^T(0)Px(0) = \sum_{k=0}^{\infty}\left[x^T(k)Qx(k) + u^T(k)Ru(k)\right] \quad (40)$$

For a PID controller as in our case, initial values of the state variables (i.e. error, its rate and integral) can not be calculated directly to find out the optimal control cost (23) since with a step-input excitation the initial value of the error rate will tend to infinity and initial value of integral error will tend to zero with the initial value of error signal remaining one. To overcome this problem the following methodology has been adopted. Here, eigen-values of the GA based differential $P$ matrices (with DARE) are evaluated, corresponding to the gains in Table II. For the oscillatory system with $\xi = 0.2$,

$$eig(P_{\Lambda=1.0} - P_{\Lambda=0.5}) = \begin{bmatrix} 1.9362 & 2.9975 & 41.0457 \end{bmatrix}^T$$
$$eig(P_{\Lambda=1.5} - P_{\Lambda=1.0}) = \begin{bmatrix} 0.7795 & 0.9049 & 11.1710 \end{bmatrix}^T \quad (41)$$

Similarly, for the sluggish system with $\xi = 5$,

$$eig(P_{\Lambda=0.5} - P_{\Lambda=1.0}) = \begin{bmatrix} 6 & 88.1 & 7128.4 \end{bmatrix}^T$$
$$eig(P_{\Lambda=1.0} - P_{\Lambda=1.5}) = \begin{bmatrix} 0.0601 & 2.6516 & 213.4166 \end{bmatrix}^T \quad (42)$$

Clearly, the eigen-values in (41) and (42) are positive which indicates that the differential matrices in left hand side of (41) and (42) are positive definite. Now, it is well known that for any two Riccati solutions $P_1, P_2$ with $P_1 > P_2$, considering initial value of the state variables as $x(0)$, pre-multiplication with $x^T(0)$ and post-multiplication with $x(0)$ yields the comparison of the associated LQR costs:

$$x^T(0)P_1x(0) > x^T(0)P_2x(0) \quad \Rightarrow J_1 > J_2 \quad (43)$$

Here, the comparison of cost of control for LQR based PID controller design has been adopted from other fractional order approaches of guaranteed dominant pole placement [35]. From (41)-(43) it is clear that while designing an optimum discrete time LQR based digital PID controller to control an oscillatory system, the cost of control increases for higher values of the fractional order integral performance index whereas, the converse is true for discrete time LQR-PID control of a sluggish system. It can also be seen that the fractional order is essentially an additional degree of freedom in the design process and for both the cases the fractional order performance index based approach is better than the integer order case. From equation (41), we find that for the cost of control $x^T(0)P_{\Lambda=1.5}x(0) > x^T(0)P_{\Lambda=1.0}x(0) > x^T(0)P_{\Lambda=0.5}x(0)$ and from equation (42), we have $x^T(0)P_{\Lambda=1.5}x(0) < x^T(0)P_{\Lambda=1.0}x(0) < x^T(0)P_{\Lambda=0.5}x(0)$. Thus it can be observed that in both cases the integer order case lies in between the fractional order cases. For the sluggish processes, $\Lambda=1.5$ gives the best design (lowest cost of control) and for oscillatory processes $\Lambda=0.5$, gives the best design. Therefore, it is recommended to set low values of the FO ($\Lambda<1$) for controlling oscillatory processes and high values of FO ($\Lambda>1$) for the control of sluggish processes in the arbitrary order cost function (31).

**5. Conclusion**

GA based optimum selection of weighting matrices have been done for the design of a discrete time LQR. Conventional PID controller design has been generalized as a LQR problem for the control of second order sluggish and oscillatory systems. The cost of control for the LQR-PID design has been shown to be dependent on the process characteristics, also the damping and the fractional order of the cost function. It is shown that the effect of inherent long memory behavior in the fractional differ-integrals can give interesting results which can be effectively harnessed to give better controller designs. Recommendation for choice of fractional order of cost function in LQR weight selection shows that the problem can be applied for the control of wide variety of industrial processes. Future scope of research can be directed towards extension of the proposed methodology to Linear Quadratic Gaussian (LQG) problems considering noisy measurements and disturbances.